\documentclass[12pt]{amsart}
\usepackage{amsmath}
\usepackage{amssymb}
\usepackage{amsmath,amssymb,amsthm,amscd}

\usepackage{txfonts}
\usepackage{amsbsy}
\usepackage{graphicx,color}

\numberwithin{equation}{section}

\newtheorem{theo}{Theorem}[section]
\newtheorem{lemm}{Lemma}[section]
\newtheorem{coro}{Corollary}[section]

\def\begeq{\begin{equation}}
\def\endeq{\end{equation}}

\begin{document}

\title{The Liouville Theorem of a Torsion System and Its Application to Symmetry Group of a Porous Medium type Equation on Symmetric Spaces}
\author{Xiao-Peng Chen \& Shi-Zhong Du \& Tian-Pei Guo}
\thanks{Author SZ is partially supported by ``Special Funds for Scientific and Technological Innovation Strategy in Guangdong Province" and  NSFG (2019A1515010605), Author XP is partially supported by the NSFC (11501344) and NSFG (2017A030313005)}
  \address{Xiao-Peng Chen, The Department of Mathematics,
            Shantou University, Shantou, 515063, P. R. China.} \email{xpchen@stu.edu.cn}
  \address{Shi-Zhong Du* (corresponding author), The Department of Mathematics,
            Shantou University, Shantou, 515063, P. R. China.} \email{szdu@stu.edu.cn}
   \address{Tian-Pei Guo, The Department of Mathematics,
            Shantou University, Shantou, 515063, P. R. China.} \email{18tpguo@stu.edu.cn}

\renewcommand{\subjclassname}{%
  \textup{2010} Mathematics Subject Classification}
\subjclass[2010]{53C35 $\cdot$ 35K59 $\cdot$ 35K65}
\date{May. 2020}
\keywords{Porous Medium Equation, Prolongation Formula}

\begin{abstract}
  In this paper, we will first prove a Liouville theorem to the torsion system
    $$
     \begin{cases}
       \xi^i_i=\lambda(x)\pm\frac{2x^k\xi^k}{|x|^2+1}, & \forall i=1,2,\cdots,n\\
       \xi^i_j+\xi^j_i=0, & \forall i\not=j
     \end{cases}
    $$
  for $(\xi,\lambda)\in C^\infty({\mathbb{R}}^n,{\mathbb{R}}^n\times{\mathbb{R}})$. As an application, complete resolutions of symmetry groups to the porous medium equation
    $$
      u_t-\triangle_g(u^m)=u^p, \ \ \forall(x,t)\in M\times{\mathbb{R}}
    $$
  of Fujita type are obtained, where $M$ is the sphere ${\mathbb{S}}^n\subset{\mathbb{R}}^{n+1}$ or hyperbolic space ${\mathbb{H}}^n$ with canonical metric $g$.
\end{abstract}

\maketitle\markboth{Porous Medium Equation}{Prolongation Formula}

\tableofcontents

\section{Introduction}

 A classical method to find  symmetry reductions of pdes is the Lie group method~\cite{O}. The classical symmetries of the general nonlinear heat equation
was considered by Clarkson and Mansfield~\cite{CM}. They gave a catalogue of
symmetry reductions.
The classical Lie group method has been  generalized  to the  porous medium
equation~\cite{V} . Gandarias~\cite{G2} applied the Lie-group  formalism to deduce symmetries of the  porous medium
equation. Franco~\cite{F}  globalize the
symmetry group of  the n-dimensional nonlinear porous medium
equation. Some other results concerning the symmetry groups of \eqref{e1.1} may refer to \cite{G,G3,G4,Q} for one spatial dimensional case, and refer to \cite{BAA,CQ,DR,GRRB} for higher dimensional case.

 In this paper, we will firstly study a system of differential equation
    \begin{equation}\label{e1.1}
     \begin{cases}
       \xi^i_i=\lambda(x)\pm\frac{2x^k\xi^k}{|x|^2+1}, & \forall i=1,2,\cdots,n\\
       \xi^i_j+\xi^j_i=0, & \forall i\not=j
     \end{cases}
    \end{equation}
for
   $$
     (\xi,\lambda)=((\xi^1,\cdots,\xi^n,\lambda)\in C^\infty({\mathbb{R}}^n,{\mathbb{R}}^n\times{\mathbb{R}}),
   $$
and then prove the following Liouville property:

\begin{theo}\label{t1.1}
 The unique solution of \eqref{e1.1} is given by $\lambda(x)\equiv0$ and
    \begin{equation}\label{e1.2}
      \xi^i=\Sigma_{j\not=i}a^i_jx^j+b^i, \ \ a^i_j=-a^j_i, \ \ \forall i\not=j,
    \end{equation}
where $a^i_j, i,j=1,2,\cdots,n$ and $b^i, i==1,2,\cdots,n$ are both constants.
\end{theo}

The system \eqref{e1.1} comes frequently from the studying of Lie's theory to PDEs on symmetric spaces. As an application, we derive classification results of symmetry groups to a porous medium equation
  \begin{equation}\label{e1.3}
    u_t-\triangle_g u^m=u^{p}, \ \ \forall(x,t)\in M\times{\mathbb{R}}
  \end{equation}
 of Fujita type, where $n$ is assumed to be greater than one and $(M,g)$ is a complete Riemannian manifold equipped with metric $g$.

  The classical frame of porous medium type equation was set up after landmark works of \cite{FdPV,LS,NST,S,S2,V}.  Although the applications of Lie's theory and the analyses of qualitative behaviors of solutions on Riemannian manifolds \cite{GM,GMP,GMP2,GMV,GMV2,V} are known for porous medium type equation, the interaction of both roles is not clear yet. It is the main purpose of this paper to clarify the symmetry group for equation on Riemmanian manifold.

  If one sets $v(x,t)=u^m(x,t/m)$ and $r=(1-m)/m, q=p/m$, \eqref{e1.3} changes to a semilinear equations in form of
   \begin{equation}\label{e1.4}
     u_t=u^{-r}(\triangle_g u+u^q).
   \end{equation}
For the sake of simplicity, we turn to classify the symmetry groups of \eqref{e1.4} with constants $q\not=0,r\not=-1$ instead of \eqref{e1.3}. In case $M={\mathbb{S}}^n\subset{\mathbb{R}}^{n+1}$, we will prove the following result.

 \begin{theo}\label{t1.2}
 When $M={\mathbb{S}}^n\subset{\mathbb{R}}^{n+1}$, after rotation on circle and the stereo polar projection from north polar, the symmetry groups of \eqref{e1.2} are generated

(1) by \eqref{e4.9}, in case $q,r+1,1$ don't equal mutually,

(2) by \eqref{e4.16}-\eqref{e4.17}, in case of $q=r+1\not=1$,

(3) by \eqref{e4.23}, in case of $q=1\not=r+1$,

(4) by \eqref{e4.29}, in case of $r+1=1\not=q$,

(5) by \eqref{e4.39}, in case of $q=r+1=1$.
\end{theo}

\noindent And in case of $M={\mathbb{R}}^n$, the following characterization result was shown.

\begin{theo}\label{t1.3}
When $M={\mathbb{H}}^n$, the symmetry groups of \eqref{e1.2} are generated by

(1) by \eqref{e5.9}, in case $q,r+1,1$ don't equal mutually,

(2) by \eqref{e5.16}-\eqref{e5.17}, in case of $q=r+1\not=1$,

(3) by \eqref{e5.23}, in case of $q=1\not=r+1$,

(4) by \eqref{e5.29}, in case of $r+1=1\not=q$,

(5) by \eqref{e5.39}, in case of $q=r+1=1$.
\end{theo}

Since there is no asymptotic assumption was imposed, our classification results Theorem \ref{t1.2}-\ref{t1.3} are complete. The contents of this paper are organized as follows. In Section 2, we recall some basic fact about the Lie's theory for PDEs. Next, the Liouville property will be proven for torsion system \eqref{e1.1} in Section 3. Finally, as an application of Theorem \ref{t1.1}, we give the proofs of Theorem \ref{t1.2} in Section 4 and Theorem \ref{t1.3} in Section 5.

\vspace{40pt}

\section{Preliminary facts to Lie's theorem on manifolds}

Let's first recall some facts of Lie's theorem (see for example \cite{O}) to partial differential equation. Here, we consider a parabolic partial differential equation
   \begin{equation}\label{e2.1}
     F(x, t, u, Du, D^2u, u_t)=0
   \end{equation}
of second order, and suppose that
  $$
   \overrightarrow{v}=\xi^i(x,t,u)\frac{\partial}{\partial x^i}+\eta(x,t,u)\frac{\partial}{\partial t}+\phi(x,t,u)\frac{\partial}{\partial u}
  $$
is an infinitesimal generator of one-parameter group action $g(\varepsilon), \varepsilon\in{\mathbb{R}}$.  By the prolongation formula in \cite{O} (Theorem 2.36, Page 110)
   \begin{eqnarray}\nonumber\label{e2.2}
     pr^{(2)}\overrightarrow{v}&=&\xi^i\frac{\partial}{\partial x^i}+\eta\frac{\partial}{\partial t}+\phi\frac{\partial}{\partial u}+\phi^i\frac{\partial}{\partial u_i}+\phi^t\frac{\partial}{\partial u_t}\\
     &&+\phi^{ij}\frac{\partial}{\partial u_{ij}}+\phi^{it}\frac{\partial}{\partial u_{it}}+\phi^{tt}\frac{\partial}{\partial u_{tt}},
   \end{eqnarray}
where
  \begin{eqnarray*}
    \phi^i&=&D_i(\phi-\xi^ju_j-\eta u_t)+\xi^ju_{ij}+\eta u_{it}\\
     &=&\phi_i+\phi_uu_i-(\xi^j_i+\xi^j_uu_i)u_j-(\eta_i+\eta_uu_i)u_t,\\
     \phi^t&=&D_t(\phi-\xi^ju_j-\eta u_t)+\xi^ju_{jt}+\eta u_{tt}\\
      &=&\phi_t+\phi_uu_t-(\xi^j_t+\xi^j_uu_t)u_j-(\eta_t+\eta_uu_t)u_t
  \end{eqnarray*}
and
   \begin{eqnarray*}
     \phi^{ij}&=&D_{ij}(\phi-\xi^ku_k-\eta u_t)+\xi^ku_{ijk}+\eta u_{ijt}\\
      &=&D_j\Bigg\{\phi_i+\phi_uu_i-(\xi^k_i+\xi^k_uu_i)u_k-\xi^ku_{ik}-(\eta_i+\eta_uu_i)u_t-\eta u_{it}\Bigg\}+\xi^ku_{ijk}+\eta u_{ijt}\\
      &=&\phi_{ij}+\phi_{iu}u_j+(\phi_{uj}+\phi_{uu}u_j)u_i+\phi_uu_{ij}-\Bigg\{\xi^k_{ij}+\xi^k_{iu}u_j+(\xi^k_{uj}+\xi^k_{uu}u_j)u_i\Bigg\}u_k\\
      &&-(\xi^k_i+\xi^k_uu_i)u_{jk}-(\xi^k_j+\xi^k_uu_j)u_{ik}-\Bigg\{\eta_{ij}+\eta_{iu}u_j+(\eta_{uj}+\eta_{uu}u_j)u_i\Bigg\}u_t\\
      &&-(\eta_i+\eta_uu_i)u_{jt}-(\eta_j+\eta_uu_j)u_{it}.
   \end{eqnarray*}
Moreover, $g(\cdot)$ is a one-parameter symmetry group of \eqref{e2.1} if and only if
   \begin{equation}\label{e2.3}
      pr^{(2)}\overrightarrow{v}F(x,t,u,Du,D^2u,u_t)=0
   \end{equation}
holds for any $u^{(2)}\equiv(u, Du, D^2u, u_t)$ satisfying
   \begin{equation}\label{e2.4}
     F(x, t, u^{(2)})=0,
   \end{equation}
where $u, Du, D^2u, u_t$ are regarded as independent variables as usually.

 When consider solution $u$ of \eqref{e1.4} on Riemannian manifold $(M,g)$ which can be parametrized by global coordinates $x\in{\mathbb{R}}^n$, we have
   \begin{equation}\label{e2.5}
     F(x,t,u,Du,D^2u,u_t)\equiv u_t-u^{-r}\Big\{g^{ij}\Big(D_{ij}u-\Gamma^k_{ij}D_ku\Big)+u^q\Big\},
   \end{equation}
 where $\Gamma^k_{ij}, i,j,k=1,2,\cdots,n$ are the Christofell symbols of $M$.

\vspace{40pt}

\section{Proof of Liouville theorem to the torsion system \eqref{e1.1}}

  In this section, we will prove that the Theorem \ref{t1.1} holds true for torsion system \eqref{e1.1}. Without some specification, we may only give the proof to
    $$
     \xi^i_i=\lambda(x)+\frac{2x^k\xi^k}{|x|^2+1},
    $$
  the minus case is similarly.

   At first, multiplying the second identity of \eqref{e1.1} by $x^ix^j$, one concludes that
   \begin{equation}\label{e3.1}
     x^jD_j(x^i\xi^i)+x^iD_i(x^j\xi^j)=0, \ \ \forall i\not=j.
   \end{equation}
Another hand, one also has that
   \begin{eqnarray}\nonumber\label{e3.2}
     x^iD_i(x^i\xi^i)&=&|x^i|^2D_i\xi^i-x^i\xi^i\\
      &=&\lambda|x^i|^2+\frac{2x^k\xi^k}{|x|^2+1}|x^i|^2-x^i\xi^i, \ \ \forall i
   \end{eqnarray}
by first formula in \eqref{e1.1}. Adding \eqref{e3.1} by two times of \eqref{e3.2} and then summing over all indices $i,j=1,2,\cdots$, one obtains that
  \begin{equation}\label{e3.3}
   x\cdot D(x\cdot\xi)=\lambda|x|^2+\frac{2(x\cdot\xi)}{|x|^2+1}|x|^2-(x\cdot\xi).
  \end{equation}
Using polar coordinates
   $$
     (r,\theta), r\equiv|x|\geq0, \theta\equiv\frac{x}{|x|}\in{\mathbb{S}}^{n-1}
   $$
and expressing $x\cdot\xi=\varphi(r,\theta)$, we have
   \begin{equation}\label{e3.4}
     \begin{cases}
       r\varphi_r=\lambda(r,\theta) r^2+\frac{r^2-1}{r^2+1}\varphi, & \forall r\geq0\\
       \varphi(0,\theta)=\varphi_0(\theta), \ \ \varphi_r(0,\theta)=0, & \forall \theta\in{\mathbb{S}}^{n-1}.
     \end{cases}
   \end{equation}
Solving the above first order O.D.E., one concludes that
   \begin{equation}\label{e3.5}
     \varphi(r,\theta)=\frac{r^2+1}{r}\Big[\int_0^r\frac{s^2}{s^2+1}\lambda(s,\theta)ds+C(\theta)\Big].
   \end{equation}
Using the initial condition of $\varphi$, there must be $C(\theta)=0, \ \forall \theta\in{\mathbb{S}}^{n-1}$. Thus,
   \begin{equation}\label{e3.6}
     x\cdot\xi(x)=\varphi(r,\theta)\equiv\frac{|x|^2+1}{|x|}\int^{|x|}_0\frac{s^2}{s^2+1}\lambda(s,\theta)ds, \forall x\in{\mathbb{R}}^n.
   \end{equation}
Next, summing \eqref{e3.1} for all $j\not=i$ and then adding it by \eqref{e3.2}, we obtain that
  $$
   x\cdot D(x^i\xi^i)+x^iD_i(x\cdot\xi)=2\lambda|x^i|^2+\frac{4x^k\xi^k}{|x|^2+1}|x^i|^2-2x^i\xi^i,
  $$
or equivalent
   \begin{eqnarray}\nonumber\label{e3.7}
    x\cdot D(x^i\xi^i)+2x^i\xi^i&=&-x^iD_i\Bigg\{\frac{r^2+1}{r}\int^r_0\frac{s^2}{s^2+1}\lambda(s,\theta)ds\Bigg\}\\
     &&+2\lambda|x^i|^2+\frac{4|x^i|^2}{r}\int^r_0\frac{s^2}{s^2+1}\lambda(s,\theta)ds\\ \nonumber
     &=&\Bigg(\frac{3}{r}+\frac{1}{r^3}\Bigg)|x^i|^2\int^r_0\frac{s^2}{s^2+1}\lambda(s,\theta)ds+|x^i|^2\lambda(r,\theta).
   \end{eqnarray}
Setting
  $$
    \psi(r,\theta)\equiv x^i\xi^i, \ \ r\equiv|x|, \ \ \theta\equiv\frac{x}{|x|},
  $$
one gets
  \begin{equation}\label{e3.8}
    r\psi_r+2\psi=\Bigg(3r+\frac{1}{r}\Bigg)\cos^2\theta_i\int^r_0\frac{s^2}{s^2+1}\lambda(s,\theta)ds+r^2\cos^2\theta_i\lambda(r,\theta).
  \end{equation}
Solving this O.D.E., we obtain that
   \begin{eqnarray}\nonumber\label{e3.9}
     & x^i\xi^i=\psi(r,\theta)=\Bigg(r+\frac{1}{r}\Bigg)\cos^2\theta_i\int^r_0\frac{s^2}{s^2+1}\lambda(s,\theta)ds&\\
     &\Leftrightarrow \xi^i=\frac{r^2+1}{r^3}x_i\int^r_0\frac{s^2}{s^2+1}\lambda(s,\theta)ds.&
   \end{eqnarray}
As a result,
  \begin{equation}\label{e3.10}
    x^jD_j(x^i\xi^i)=\frac{x_i^2x_j^2}{r^3}\Bigg\{-\Bigg(1+\frac{3}{r^2}\Bigg)\int^r_0\frac{s^2}{s^2+1}\lambda(s,\theta)ds+r\lambda(r,\theta)\Bigg\}.
  \end{equation}
Substituting into \eqref{e3.1}, one concludes that
   \begin{equation}\label{e3.11}
     \int^r_0\frac{s^2}{s^2+1}\lambda(s,\theta)ds=\frac{r^3}{r^2+3}\lambda(r,\theta).
   \end{equation}
Therefore, $\lambda$ is the solution to
  \begin{equation}\label{e3.12}
    \begin{cases}
      \partial_r\lambda=-\frac{4r}{(r^2+1)(r^2+3)}\lambda, & \forall r>0\\
      \lambda(0,\theta)=\lambda_0, \ \ \partial_r\lambda(0,\theta)=0,
    \end{cases}
  \end{equation}
which is given
   \begin{equation}\label{e3.13}
    \lambda(r,\theta)=\lambda_0\frac{r^2+3}{r^2+1}.
   \end{equation}
Taking derivative on second identity of \eqref{e3.9}, and comparing both sides of \eqref{e1.1}, we derive that
  \begin{equation}\label{e3.14}
    \frac{1}{r^3}\Bigg\{1-r^2-x_ir-\frac{3x_i}{r}\Bigg\}\int^r_0\frac{s^2}{s^2+1}\lambda(s,\theta)ds=\frac{r-x_i}{r}\lambda(r,\theta)
  \end{equation}
for function given by \eqref{e3.13}. So, there must be $\lambda_0=0$ and thus $\lambda(x)\equiv0, \ \forall x\in{\mathbb{R}}^n$. Now, the conclusion of Theorem \ref{t1.1} follows from the following Lemma.

 \begin{lemm}\label{l3.1}
  The unique solutions $(\lambda,\xi)\in{\mathbb{R}}\times C^\infty({\mathbb{R}}^n,{\mathbb{R}}^n)$ of
    \begin{equation}\label{e3.15}
     \begin{cases}
       \xi^i_i=\lambda+\frac{2x^k\xi^k}{|x|^2+1}, & \forall i\\
       \xi^i_j+\xi^j_i=0, & \forall i\not=j
     \end{cases}
    \end{equation}
  are given by $\lambda=0$ and
    \begin{equation}\label{e3.16}
      \xi^i=\Sigma_{j\not=i}a^i_jx^j+b^i, \ \ a^i_j=-a^j_i, \ \ \forall i\not=j.
    \end{equation}
\end{lemm}

\noindent\textbf{Remark.} Taking derivatives on \eqref{e3.15}, one gets
    \begin{equation}\label{e3.17}
     \begin{cases}
     \xi^j_{ii}=-2D_j\Big(\frac{x^k\xi^k}{|x|^2+1}\Big), & \forall i\not=j\\
     \xi^j_{jj}=2D_j\Big(\frac{x^k\xi^k}{|x|^2+1}\Big), & \forall j
     \end{cases}
    \end{equation}
and hence
    \begin{equation}\label{e3.18}
      \triangle\xi^j+(n-2)D_j^2\xi^j=0, \ \ \forall j.
    \end{equation}
However, since there is no asymptotic assumption of $\xi$ at infinity known, it would be difficulty to classify all solutions $\xi$ using harmonicity formula \eqref{e3.18}. Fortunately, using the special structure of system of differential equations \eqref{e3.15}, one can prove the Lemma \ref{l3.1}.\\

\noindent\textbf{Proof.} Multiplying the second identity of \eqref{e3.15} by $x^ix^j$, one gets
   \begin{equation}\label{e3.19}
     x^jD_j(x^i\xi^i)+x^iD_i(x^j\xi^j)=0, \ \ \forall i\not=j.
   \end{equation}
Another hand, we have
   \begin{eqnarray}\nonumber\label{e3.20}
     x^iD_i(x^i\xi^i)&=&|x^i|^2D_i\xi^i-x^i\xi^i\\
      &=&\lambda|x^i|^2+\frac{2x^k\xi^k}{|x|^2+1}|x^i|^2-x^i\xi^i, \ \ \forall i
   \end{eqnarray}
by first formula in \eqref{e3.15}. Adding \eqref{e3.17} by two times of \eqref{e3.18} and then summing over all indices $i,j$, one obtains that
  \begin{equation}\label{e3.21}
   x\cdot D(x\cdot\xi)=\lambda|x|^2+\frac{2(x\cdot\xi)}{|x|^2+1}|x|^2-(x\cdot\xi).
  \end{equation}
Using polar coordinates
   $$
     (r,\theta), r\equiv|x|\geq0, \theta\equiv\frac{x}{|x|}\in{\mathbb{S}}^{n-1}
   $$
and expressing $x\cdot\xi=\varphi(r,\theta)$, we have
   \begin{equation}\label{e3.22}
     \begin{cases}
       r\varphi_r=\lambda r^2+\frac{r^2-1}{r^2+1}\varphi, & \forall r\geq0\\
       \varphi(0,\theta)=\varphi_0(\theta), \ \ \varphi_r(0,\theta)=0, & \forall \theta\in{\mathbb{S}}^{n-1}.
     \end{cases}
   \end{equation}
Solving the above first order O.D.E., one concludes that
   \begin{equation}\label{e3.23}
     \varphi(r,\theta)=\frac{r^2+1}{r}\Big[\lambda\Big(r-\arctan r\Big)+C(\theta)\Big].
   \end{equation}
Using the initial condition of $\varphi$, there must be $C(\theta)=0, \ \forall \theta\in{\mathbb{S}}^{n-1}$. Thus,
   \begin{equation}\label{e3.24}
     x\cdot\xi(x)=\varphi(r,\theta)\equiv\lambda\frac{|x|^2+1}{|x|}(|x|-\arctan|x|), \forall x\in{\mathbb{R}}^n.
   \end{equation}
It's clear that for solution $\xi$ of \eqref{e3.15}, \eqref{e3.17} and \eqref{e3.18} also hold true. Therefore, one concludes that
   \begin{eqnarray}\nonumber\label{e3.25}
     &&\Bigg(\frac{\partial^2}{\partial (x^i)^2}+\frac{\partial^2}{\partial (x^j)^2}\Bigg)D_j\Bigg(\frac{x\cdot\xi}{|x|^2+1}\Bigg)=0\\
     &\Leftrightarrow&\Bigg(\frac{\partial^2}{\partial (x^i)^2}+\frac{\partial^2}{\partial (x^j)^2}\Bigg)\Bigg\{-\lambda\Bigg(\frac{1}{|x|(|x|^2+1)}-\frac{\arctan|x|}{|x|^2}\Bigg)\frac{x^j}{|x|}\Bigg\}=0.
   \end{eqnarray}
Using Taylor's expansion for function
   $$
    f(r)\equiv-\frac{1}{r}\Bigg\{\frac{1}{r(r^2+1)}-\frac{\arctan r}{r^2}\Bigg\}=\Sigma_{k=0}^\infty(-1)^{k+1}\frac{2k+2}{2k+3}r^{2k},
   $$
we know that
   $$
    F(x^i,x^j)\equiv-\lambda\Bigg(\frac{1}{|x|(|x|^2+1)}-\frac{\arctan|x|}{|x|^2}\Bigg)\frac{x^j}{|x|}\in C^\infty({\mathbb{R}}^2)
   $$
is a smooth harmonic function on $(x^i,x^j)\in{\mathbb{R}}^2$. Noting that it is also bounded, a famous Liouville theorem \cite{GT} for harmonic function gives the unique possibility by $\lambda=0$ and thus
  \begin{equation}\label{e3.26}
    x\cdot\xi(x)=0, \ \ \forall x\in{\mathbb{R}}^n.
  \end{equation}
Combining with \eqref{e3.17}, it yields that
    $$
     \xi^i=\Sigma_{j}a^i_jx^j+b^i
    $$
and
    $$
     a^i_j=-a^j_i, \ \ a^i_i=0, \ \ \forall i\not=j.
    $$
The proof of the lemma was done. $\Box$\\

Thus, we have completed the proof of Theorem \ref{t1.1}. $\Box$\\

\vspace{40pt}

\section{Symmetry Group on ${\mathbb{S}}^n$}

Let's start with a well known Hairy ball theorem (refer to a analytic proof by Milnor \cite{M}).

\begin{theo}\label{t3.1}
  Suppose that $n$ is even, then there is no non-vanishing tangential vector field on ${\mathbb{S}}^n$.
\end{theo}

As a corollary, we have the following Brouwer type fix point theorem on sphere.

\begin{coro}\label{c3.1}
  Suppose that $n$ is even, then for any continuous mapping $f: {\mathbb{S}}^n\to{\mathbb{S}}^n$ satisfying
     $$
      f(x)\not=-x, \ \ \forall x\in{\mathbb{S}}^n,
     $$
  there is at least one fix point $x_0\in{\mathbb{S}}^n$.
\end{coro}

\noindent\textbf{Proof.} Suppose on the contrary, there is a continuous mapping $f: {\mathbb{S}}^n\to{\mathbb{S}}^n$ satisfying
   $$
    f(x)\not=x, -x, \ \ \forall x\in{\mathbb{S}}^n.
   $$
Setting
   $$
    v(x)\equiv f(x)-(f(x),x)x, \ \ \forall x\in{\mathbb{S}}^n,
   $$
we have $v$ is a tangential vector field on ${\mathbb{S}}^n$. Moreover, $v$ does not vanish everywhere on ${\mathbb{S}}^n$ due to $|(f(x),x)|<1$ for all $x\in{\mathbb{S}}^n$, which contradicts with the hairy ball theorem \ref{t3.1}. The proof was done. $\Box$\\

\noindent\textbf{Remark.} The assumption $f(x)\not=-x, \ \forall x\in{\mathbb{S}}^n$ can not be removed, since the inversion transformation is a continuous mapping of ${\mathbb{S}}^n$ which has no fix point.\\

By Corollary \ref{c3.1}, we may assume the symmetry group action has at least one fix point for $n$ is even, as well as for $n$ is odd after a rotation $SO(1)$. Without loss of generality, we may also assume that the fix point is $(0,1)\in{\mathbb{R}}^{n}\times{\mathbb{R}}$. Letting $X=(y,z)\in{\mathbb{R}}^n\times{\mathbb{R}}$ be the coordinate representation of ${\mathbb{S}}^n$, we have the spherical polar projection $P: {\mathbb{S}}^n\to{\mathbb{R}}^n$ is given by
   $$
    P(y,z)\equiv x=\frac{y}{1-z}, \ \ \forall (y,z)\in {\mathbb{S}}^n\setminus(0,1).
   $$
Regarding $x\in{\mathbb{R}}^n$ as local coordinates of ${\mathbb{S}}^n\setminus(0,1)$, one can express $X$ in terms of $x$ by
   $$
    X=\Bigg(\frac{2x}{|x|^2+1},\frac{|x|^2-1}{|x|^2+1}\Bigg).
   $$
Therefore, upon this local coordinates, the induced metric $g$ of ${\mathbb{S}}^n$ is given by
   $$
    g_{ij}\equiv\Bigg(\frac{\partial X}{\partial x^i},\frac{\partial X}{\partial x^j}\Bigg)=\frac{4\delta_{ij}}{(|x|^2+1)^2}
   $$
and the Chritoffel symbol
   \begin{eqnarray*}
    \Gamma^k_{ij}&=&\frac{1}{2}g^{kl}\Bigg(\frac{\partial g_{jl}}{\partial x^i}+\frac{\partial g_{il}}{\partial x^j}-\frac{\partial g_{ij}}{\partial x^l}\Bigg)\\
     &=&\frac{2}{|x|^2+1}\Big(-x_i\delta_{jk}-x_j\delta_{ik}+x_k\delta_{ij}\Big).
   \end{eqnarray*}
Thus, the Laplace-Beltrami operator
   \begin{eqnarray*}
    \triangle_g&=&\frac{1}{\sqrt{g}}\frac{\partial}{\partial x^i}\Bigg(\sqrt{g}g^{ij}\frac{\partial}{\partial x^j}\Bigg)=g^{ij}\Bigg(\frac{\partial^2}{\partial x^i\partial x^j}-\Gamma^k_{ij}\frac{\partial}{\partial x^k}\Bigg)\\
    &=&\frac{(|x|^2+1)^2}{4}\triangle-(n-2)\frac{|x|^2+1}{2}x\cdot\nabla,
   \end{eqnarray*}
where
   $$
     \triangle\equiv\Sigma_i\frac{\partial^2}{\partial (x^i)^2}, \ \ \nabla\equiv\Bigg(\frac{\partial}{\partial x^1},\cdots,\frac{\partial}{\partial x^n}\Bigg)
   $$
are the canonical Laplace operator and gradient operator of ${\mathbb{R}}^n$  under flat metric respectively.

Suppose that
  $$
   \overrightarrow{v}=\xi^i(x,t,u)\frac{\partial}{\partial x^i}+\eta(x,t,u)\frac{\partial}{\partial t}+\phi(x,t,u)\frac{\partial}{\partial u}
  $$
to be an infinitesimal generator of one-parameter symmetry group $g(\varepsilon), \varepsilon\in{\mathbb{R}}$, we have the prolongation formula is given by \eqref{e2.2}. Using spherical polar projection, \eqref{e2.5} changes to
  \begin{equation}\label{e4.1}
    u^ru_t=\frac{(|x|^2+1)^2}{4}\triangle u-(n-2)\frac{|x|^2+1}{2}x\cdot\nabla u+u^q.
  \end{equation}
Therefore, the group action $g(\cdot)$ is a one-parameter symmetry group of \eqref{e4.1} if and only if
   \begin{eqnarray}\nonumber\label{e4.2}
     &ru^{r-1}u_t\phi+u^r\phi^t=(|x|^2+1)x^i\xi^i\triangle u+\frac{(|x|^2+1)^2}{4}\phi^{ii}&\\
     &-(n-2)\Bigg\{x^i\xi^ix\cdot\nabla u+\frac{|x|^2+1}{2}u_i\xi^i+\frac{|x|^2+1}{2}x^i\phi^i\Bigg\}+qu^{q-1}\phi.&
   \end{eqnarray}
Comparing the like terms on both sides of \eqref{e4.2}, one gets that
  \begin{eqnarray*}
    &\eta=\eta(t)&\\
    &\xi^k_u=0, \ \ \forall k&\\
    &\phi_{uu}=0&\\
    &ru^{-1}\phi-\eta_t-\frac{4x^k}{|x|^2+1}\xi^k+2\xi^i_i=0, \ \ \forall i&\\
    &\xi^i_j+\xi^j_i=0, \ \ \forall i\not=j&\\
    &\frac{2}{|x|^2+1}u^r\xi^k_t+\frac{|x|^2+1}{2}(2\phi_{ku}-\triangle\xi^k)-ru^{-1}\phi x^k&\\ &+\Big(\eta_t+(1-n)\phi_u\Big)x^k+(2-n)\Bigg\{\frac{2x^i}{|x|^2+1}\xi^ix^k+\xi^k-\xi^k_ix^i\Bigg\}=0, \ \ \forall k&\\
    &(r-q)u^{q-1}\phi+u^r\phi_t+(\phi_u-\eta_t)u^q-\frac{(|x|^2+1)^2}{4}\triangle\phi-(2-n)\frac{|x|^2+1}{2}x^i\phi_i=0,&
  \end{eqnarray*}
or equivalent to
   \begin{eqnarray}\nonumber\label{e4.3}
     &\xi^k=\xi^k(x,t),\ \ \eta=\eta(t), \ \ \phi=\alpha(x,t)u&\\ \nonumber
     &2\xi^i_{i}=-r\alpha(x,t)+\eta_t+\frac{4x^k\xi^k}{|x|^2+1}, \ \ \forall i&\\
     &\xi^i_j+\xi^j_i=0, \ \ \forall i\not=j&\\ \nonumber
     &\Big[(r+1-q)\alpha(x,t)-\eta_t\Big]u^q+\alpha_t(x,t)u^{r+1}-u\triangle_g\alpha=0&\\ \nonumber
     &u^r\partial_t\xi^k-\triangle_g\xi^k+\frac{(|x|^2+1)^2}{2}D_k\alpha+\frac{|x|^2+1}{2}\Big\{\eta_t+(1-n-r)\alpha\Big\}x^k&\\ \nonumber
     &+(2-n)\Bigg\{(x^i\xi^i)x^k+\frac{|x|^2+1}{2}\xi^k\Bigg\}=0, \ \ \forall k&
   \end{eqnarray}
So, our situations are divided into five cases:\\

\noindent\textbf{Case 1:} ($q, r+1, 1$ don't equal mutually)

By 4th-5th identities in \eqref{e3.3}, we get
   \begin{equation}\label{e3.4}
     \alpha=\frac{\beta}{r+1-q}, \ \ \eta(t)=\beta t+\gamma,\ \ \xi^k=\xi^k(x)
   \end{equation}
and
   \begin{equation}\label{e3.5}
     \begin{cases}
       \xi^i_i=-\frac{\beta(1-q)}{2(r+1-q)}+\frac{2x^k}{|x|^2+1}\xi^k, & \forall i\\
       \xi^i_j+\xi^j_i=0, & \forall i\not=j
     \end{cases}
   \end{equation}
and
   \begin{equation}\label{e3.6}
     \triangle_g\xi-(2-n)\frac{\xi}{1-z}=\Bigg\{\frac{\beta(2-n-q)}{r+1-q}+(2-n)y\cdot\xi\Bigg\}\frac{y}{(1-z)^2}, \ \ \xi\equiv(\xi^1,\cdots,\xi^n).
   \end{equation}
Applying Theorem \ref{t1.1} to \eqref{e4.4} and \eqref{e4.5}, one concludes that
  \begin{equation}\label{e4.7}
    \alpha=0, \ \ \eta(t)=\gamma, \ \ \xi^i=\Sigma_{j\not=i}a^i_jx^j, \ \ a^i_j=-a^j_i, \ \ \forall i\not=j.
  \end{equation}
Therefore the infinite generators are spanned by
  \begin{eqnarray}\nonumber\label{e4.8}
    \overrightarrow{v}_1&=&\partial_t\\
    \overrightarrow{v}_2&=&\Sigma_{j\not=i}a^i_jx^j\partial_{x^i},
  \end{eqnarray}
and their corresponding group actions are generated by
   \begin{eqnarray}\nonumber\label{e4.9}
     g_1(\varepsilon) &: & (x,t,u)\to(\widetilde{x},\widetilde{t},\widetilde{u})=(x,t+\varepsilon, u)\\
     g_2(\varepsilon) &: & (x,t,u)\to(\widetilde{x},\widetilde{t},\widetilde{u})=(Ax,t, u),\ \ A\in SO(n)
   \end{eqnarray}
corresponding to translation in time and rotation on space.\\

\vspace{10pt}

\noindent\textbf{Case 2:} ($q=r+1\not=1$)

By 4th identity in \eqref{e4.3}, one gets that
   \begin{eqnarray}\nonumber\label{e4.10}
     &\alpha(x,t)=\eta(t)+\beta(x), \ \ \eta(t)=\kappa e^{rt}-\gamma,\ \ \xi^k=\xi^k(x)&\\
     &\triangle_g\beta=0&
   \end{eqnarray}
and
   \begin{equation}\label{e4.11}
     \begin{cases}
       \xi^i_i=\frac{r}{2}\Big[\gamma-\beta(x)\Big]+\frac{2x^k}{|x|^2+1}\xi^k, & \forall i\\
       \xi^i_j+\xi^j_i=0, & \forall i\not=j
     \end{cases}
   \end{equation}
and
   \begin{eqnarray}\nonumber\label{e4.12}
    &\triangle_g\xi-(2-n)\frac{\xi}{1-z}=\Bigg\{(1-n)\kappa e^{rt}+(1-n-r)\Big[\beta(x)-\gamma\Big]+(2-n)y\cdot\xi\Bigg\}\frac{y}{(1-z)^2}&\\
    &+\frac{2D\beta}{(1-z)^2},\ \ D\beta\equiv(\beta_1,\cdots,\beta_n).&
   \end{eqnarray}
Applying Theorem \ref{t1.1} to \eqref{e4.11}, one concludes that $\beta=\gamma$ and
   \begin{equation}\label{e4.13}
     \xi^k=a^k_jx^j+b^k, \ \ \forall k.
   \end{equation}
So, we get
  \begin{equation}\label{e4.14}
    \alpha=\kappa e^{rt}, \ \ \eta=\kappa e^{rt}-\gamma, \ \ \xi^k=\Sigma_{j\not=k}a^k_jx^j. \ \ a^i_j=-a^j_i, \ \ \forall i\not=j.
  \end{equation}
Thus, the infinite generators are spanned by
  \begin{eqnarray}\nonumber\label{e4.15}
    \overrightarrow{v}_1&=&e^{rt}\partial_t+e^{rt}u\partial_u\\
    \overrightarrow{v}_2&=&\partial_t\\ \nonumber
    \overrightarrow{v}_3&=&\Sigma_{j\not=i}a^i_jx^j\partial_{x^i},
  \end{eqnarray}
and their corresponding group actions are generated by
  \begin{equation}\label{e4.16}
     g_1(\varepsilon): \ \left(
    \begin{array}{c}
      x\\
      t \\
      u
    \end{array}
   \right)\to \left(
    \begin{array}{c}
      \widetilde{x}\\
      \widetilde{t}\\
      \widetilde{u}
    \end{array}
   \right)=   \left(
    \begin{array}{c}
      x\\
      -\frac{1}{r}\ln(e^{-rt}-r\varepsilon)\\
      u(e^{-rt}-r\varepsilon)^{-\frac{1}{r}}e^{-t}
    \end{array}
   \right)
  \end{equation}
and
   \begin{eqnarray}\nonumber\label{e4.17}
     g_2(\varepsilon) &: & (x,t,u)\to(\widetilde{x},\widetilde{t},\widetilde{u})=(x,t+\varepsilon, u)\\
     g_3(\varepsilon) &: & (x,t,u)\to(\widetilde{x},\widetilde{t},\widetilde{u})=(A_\varepsilon x, t, u), \ \ A_\varepsilon\in SO(n)
   \end{eqnarray}

\vspace{10pt}

\noindent\textbf{Case 3:} ($q=1\not=r+1$)

 We have $\xi^k=\xi^k(x)$ and
   \begin{eqnarray}\nonumber\label{e4.18}
      &\alpha=\alpha(x),\ \ \eta=\beta t+\gamma, \ \ \xi^k=\xi^k(x)&\\
      &\triangle_g\alpha-r\alpha+\beta=0&
   \end{eqnarray}
 and
   \begin{equation}\label{e4.19}
     \begin{cases}
       2\xi^i_i=-r\alpha(x)+\beta+\frac{4x^k\xi^k}{|x|^2+1}, & \forall i\\
       \xi^i_j+\xi^j_i=0, & \forall i\not=j
     \end{cases}
   \end{equation}
 and
    \begin{equation}\label{e4.20}
     \triangle_g\xi-(2-n)\frac{\xi}{1-z}=\Big\{\beta+(2-n-r)\alpha+y\cdot\xi\Big\}\frac{y}{(1-z)^2}+\frac{2D\alpha}{(1-z)^2}.
    \end{equation}
By Theorem \ref{t1.1} and \eqref{e4.19}, one obtains that
  \begin{equation}\label{e4.21}
    \alpha=\frac{\beta}{r},\ \ \eta=\beta t+\gamma, \ \ \xi^k=\Sigma_{j\not=k}a^k_jx^j, \ \ a^i_j=-a^j_i, \ \ \forall i\not=j.
  \end{equation}
Thus, the infinite generators are spanned by
  \begin{eqnarray}\nonumber\label{e4.22}
    \overrightarrow{v}_1&=&t\partial_t+u/r\partial_u\\
    \overrightarrow{v}_2&=&\partial_t\\ \nonumber
    \overrightarrow{v}_3&=&\Sigma_{j\not=i}a^i_jx^j\partial_{x^i},
  \end{eqnarray}
and their corresponding group actions are generated by
   \begin{eqnarray}\nonumber\label{e4.23}
     g_1(\varepsilon) &: & (x,t,u)\to(\widetilde{x},\widetilde{t},\widetilde{u})=(x,e^{\varepsilon}t, e^{\varepsilon/r}u)\\
     g_2(\varepsilon) &: & (x,t,u)\to(\widetilde{x},\widetilde{t},\widetilde{u})=(x,t+\varepsilon, u)\\ \nonumber
     g_3(\varepsilon) &: & (x,t,u)\to(\widetilde{x},\widetilde{t},\widetilde{u})=(A_\varepsilon x, t, u), \ \ A_\varepsilon\in SO(n)
   \end{eqnarray}

\vspace{10pt}

\noindent\textbf{Case 4:} ($r+1=1\not=q$)

We have $\xi^k=\xi^k(x,t)$ and
  \begin{eqnarray}\nonumber\label{e4.24}
    &&\begin{cases}
      (1-q)\alpha(x,t)=\eta_t\\
      \alpha_t-\frac{(|x|^2+1)^2}{4}\triangle\alpha-\frac{|x|^2+1}{2}x\cdot\nabla\alpha=0
    \end{cases}\\
   &\Rightarrow& \eta=\beta t+\gamma, \ \ \alpha=\frac{\beta}{1-q}
  \end{eqnarray}
and
   \begin{equation}\label{e4.25}
     \begin{cases}
       2\xi^i_i=\beta+\frac{4x^k}{|x|^2+1}\xi^k, & \forall i\\
       \xi^i_j+\xi^j_i=0, & \forall i\not=j
     \end{cases}
   \end{equation}
and
  \begin{equation}\label{e4.26}
   \partial_t\xi-\triangle_g\xi+(2-n)\frac{\xi}{1-z}=-\Bigg\{\frac{\beta(2-n-q)}{1-q}+(2-n)y\cdot\xi\Bigg\}\frac{y}{(1-z)^2}.
  \end{equation}
Applying Theorem \ref{t1.1} to \eqref{e4.25}, one concludes that
  \begin{equation}\label{e4.27}
    \alpha=0, \ \ \eta=\gamma, \ \ \xi^k=\Sigma_{j\not=k}a^k_j(t)x^j,\ \ a^k_j(t)=-a^j_k(t), \ \ \forall j\not=k.
  \end{equation}
Substituting into \eqref{e4.26}, there holds
   $$
    \Sigma_{j\not=k}\partial_ta^k_j(t)x^j=0, \ \ \forall k
   $$
or equivalently $\xi^k=\xi^k(x), \forall k$. So, the infinite generators are spanned by
  \begin{eqnarray}\nonumber\label{e4.28}
    \overrightarrow{v}_1&=&\partial_t\\
    \overrightarrow{v}_2&=&\Sigma_{j\not=i}a^i_jx^j\partial_{x^i},
  \end{eqnarray}
and their corresponding group actions are generated by
   \begin{eqnarray}\nonumber\label{e4.29}
     g_1(\varepsilon) &: & (x,t,u)\to(\widetilde{x},\widetilde{t},\widetilde{u})=(x,t+\varepsilon, u)\\
     g_2(\varepsilon) &: & (x,t,u)\to(\widetilde{x},\widetilde{t},\widetilde{u})=(Ax,t, u),\ \ A\in SO(n)
   \end{eqnarray}
corresponding to translation in time and rotation on space respectively.\\

\vspace{10pt}

\noindent\textbf{Case 5:} ($q=r+1=1$)

We have
 \begin{equation}\label{e4.30}
   \alpha_t-\triangle_g\alpha=\eta_t
 \end{equation}
and
   \begin{equation}\label{e4.31}
     \begin{cases}
       2\xi^i_i=\eta_t+\frac{4x^k}{|x|^2+1}\xi^k, & \forall i\\
       \xi^i_j+\xi^j_i=0, & \forall i\not=j
     \end{cases}
   \end{equation}
and
  \begin{equation}\label{e4.32}
  \partial_t\xi-\triangle_g\xi+(2-n)\frac{\xi}{1-z}=-\Bigg\{\eta_t+(1-n)\alpha+(2-n)y\cdot\xi\Bigg\}\frac{y}{(1-z)^2}-\frac{2D\alpha}{(1-z)^2}.
  \end{equation}
Using again Theorem \ref{t1.1}, it's inferred from \eqref{e4.30}-\eqref{e4.31} that
  \begin{equation}\label{e4.33}
    \alpha_t-\triangle_g\alpha=0, \ \ \eta=\gamma,\ \ \xi^k=\Sigma_{j\not=k}a^k_j(t)x^j, \ \ a^k_j(t)=-a^j_k(t), \ \  \forall j\not=k.
  \end{equation}
Combining with \eqref{e4.32}, one concludes that
  \begin{equation}\label{e4.34}
   D_k\alpha+\frac{(1-n)x^k}{|x|^2+1}\alpha=-\Sigma_{j\not=k}\frac{(|x|^2+1)^2}{2}\partial_ta^k_j(t)x^j.
  \end{equation}
Taking a second derivative on $x^l$, there holds
  \begin{eqnarray}\nonumber\label{e4.35}
    D_{kl}\alpha&=&-(1-n)\Bigg\{\frac{\delta_{kl}}{|x|^2+1}-\frac{2x^kx^l}{(|x|^2+1)^2}\Bigg\}\alpha-(1-n)\frac{x^k}{|x|^2+1}D_l\alpha\\
     &&-2(|x|^2+1)\Sigma_{j\not=k, j\not=l}x^lx^j\partial_ta^k_j(t)-\Bigg\{\frac{(|x|^2+1)^2}{2}+2(x^l)^2(|x|^2+1)\Bigg\}\partial_ta^k_l(t).
  \end{eqnarray}
Substituting \eqref{e4.34} into \eqref{e4.35}, one derives that
   \begin{eqnarray*}
    &D_{kl}\alpha+(1-n)\Bigg\{\frac{\delta_{kl}}{|x|^2+1}+(n-3)\frac{x^kx^l}{(|x|^2+1)^2}\Bigg\}\alpha&\\
    &=\frac{1-n}{2}(|x|^2+1)\Sigma_{j\not=k, j\not=l}x^kx^j\partial_ta^l_j(t)-2(|x|^2+1)\Sigma_{j\not=k, j\not=l}x^lx^j\partial_ta^k_j(t)&\\
    &-\frac{1-n}{2}(|x|^2+1)(x^k)^2\partial_ta^k_k(t)-(|x|^2+1)\Bigg\{\frac{|x|^2+1}{2}+2(x^l)^2\Bigg\}\partial_ta^k_l(t).&
   \end{eqnarray*}
Exchanging the indices $k,l$ and using the symmetry of second derivatives, it yields that
   \begin{eqnarray*}
     &\frac{5-n}{2}\Sigma_{j\not=k, j\not=l}x^j\Bigg\{x^k\partial_ta^l_j(t)-x^l\partial_ta^k_j(t)\Bigg\}&\\
     &-\Bigg\{(|x|^2+1)+\frac{5-n}{2}\Big[(x^k)^2+(x^l)^2\Big]\Bigg\}\partial_ta^k_l(t)=0, \ \ \forall k\not=l&.
   \end{eqnarray*}
Consequently, we get
   \begin{equation}\label{e4.36}
     a^k_j(t)=a^k_j, \ \ \forall j\not=k
   \end{equation}
and so
   \begin{equation}\label{e4.37}
     \alpha=\beta, \ \ \eta=\gamma, \xi^k=\Sigma_{j\not=k}a^k_jx^j, \ \ a^k_j=-a^j_k, \ \ \forall j\not=k.
   \end{equation}
Thus, the infinite generators are spanned by
  \begin{eqnarray}\nonumber\label{e4.38}
    \overrightarrow{v}_1&=&u\partial_u\\
    \overrightarrow{v}_2&=&\partial_t\\ \nonumber
    \overrightarrow{v}_3&=&\Sigma_{j\not=i}a^i_jx^j\partial_{x^i},
  \end{eqnarray}
and their corresponding group actions are generated by
   \begin{eqnarray}\nonumber\label{e4.39}
     g_1(\varepsilon) &: & (x,t,u)\to(\widetilde{x},\widetilde{t},\widetilde{u})=(x,t,e^{\varepsilon}u)\\
     g_2(\varepsilon) &: & (x,t,u)\to(\widetilde{x},\widetilde{t},\widetilde{u})=(x,t+\varepsilon,u)\\ \nonumber
     g_3(\varepsilon) &: & (x,t,u)\to(\widetilde{x},\widetilde{t},\widetilde{u})=(A_\varepsilon x, t, u), \ \ A_\varepsilon\in SO(n)
   \end{eqnarray}

\vspace{40pt}

\section{Symmetry Group on ${\mathbb{H}}^n$}

As one knows, Hyperbolic space ${\mathbb{H}}^n, n\geq2$ is a complete, simply connected Riemannian manifolds having constant sectional curvature $-1$. The most important models of hyperbolic spaces include the Poincar\'{e} balls, half-spaces and hyperboloids (or Lorentz model). In this section, we only discuss the Poincar\'{e} balls $B_1\subset{\mathbb{R}}^n, n\geq2$, which equipped with a complete metric
   $$
    g_{ij}=\Bigg(\frac{2}{1-|x|^2}\Bigg)^2\delta_{ij}, \ \ \forall i,j=1,2,\cdots,n, \ \ \forall x\in B_1,
   $$
whose Chritoffel symbol is given by
   \begin{eqnarray*}
    \Gamma^k_{ij}&=&\frac{1}{2}g^{kl}\Bigg(\frac{\partial g_{jl}}{\partial x^i}+\frac{\partial g_{il}}{\partial x^j}-\frac{\partial g_{ij}}{\partial x^l}\Bigg)\\
     &=&\frac{2}{1-|x|^2}\Big(x^i\delta_{jk}+x^j\delta_{ik}-x^k\delta_{ij}\Big).
   \end{eqnarray*}
Thus, the Laplace-Beltrami operator
   \begin{eqnarray*}
    \triangle_g&=&\frac{1}{\sqrt{g}}\frac{\partial}{\partial x^i}\Bigg(\sqrt{g}g^{ij}\frac{\partial}{\partial x^j}\Bigg)=g^{ij}\Bigg(\frac{\partial^2}{\partial x^i\partial x^j}-\Gamma^k_{ij}\frac{\partial}{\partial x^k}\Bigg)\\
    &=&\frac{(1-|x|^2)^2}{4}\triangle+(n-2)\frac{1-|x|^2}{2}x\cdot\nabla, \ \ \forall x\in B_1
   \end{eqnarray*}
and \eqref{e2.5} changes to
  \begin{equation}\label{e5.1}
    u^ru_t=\frac{(1-|x|^2)^2}{4}\triangle u+(n-2)\frac{1-|x|^2}{2}x\cdot\nabla u+u^q.
  \end{equation}
As above, suppose that
  $$
   \overrightarrow{v}=\xi^i(x,t,u)\frac{\partial}{\partial x^i}+\eta(x,t,u)\frac{\partial}{\partial t}+\phi(x,t,u)\frac{\partial}{\partial u}
  $$
is an infinitesimal generator of one-parameter symmetry group $g(\varepsilon), \varepsilon\in{\mathbb{R}}$, the prolongation $pr^{(2)}\overrightarrow{v}$ is given by \eqref{e2.2}. Therefore, the group action $g(\cdot)$ is a one-parameter symmetry group of \eqref{e5.1} if and only if
   \begin{eqnarray}\nonumber\label{e5.2}
     &ru^{r-1}u_t\phi+u^r\phi^t=-(1-|x|^2)x^i\xi^i\triangle u+\frac{(1-|x|^2)^2}{4}\phi^{ii}&\\
     &+(n-2)\Bigg\{-x^i\xi^ix\cdot\nabla u+\frac{1-|x|^2}{2}u_i\xi^i+\frac{1-|x|^2}{2}x^i\phi^i\Bigg\}+qu^{q-1}\phi.&
   \end{eqnarray}
Comparing the like terms on both sides of \eqref{e5.2}, we get
  \begin{eqnarray*}
    &\eta=\eta(t)&\\
    &\xi^k_u=0, \ \ \forall k&\\
    &\phi_{uu}=0&\\
    &ru^{-1}\phi-\eta_t+\frac{4x^k}{1-|x|^2}\xi^k+2\xi^i_i=0, \ \ \forall i&\\
    &\xi^i_j+\xi^j_i=0, \ \ \forall i\not=j&\\
    &-\frac{2}{1-|x|^2}u^r\xi^k_t-\frac{1-|x|^2}{2}(2\phi_{ku}-\triangle\xi^k)-ru^{-1}\phi x^k&\\ &+\Big(\eta_t+(1-n)\phi_u\Big)x^k+(2-n)\Bigg\{-\frac{2x^i}{1-|x|^2}\xi^ix^k+\xi^k-\xi^k_ix^i\Bigg\}=0, \ \ \forall k&\\
    &(r-q)u^{q-1}\phi+u^r\phi_t+(\phi_u-\eta_t)u^q-\frac{(1-|x|^2)^2}{4}\triangle\phi+(2-n)\frac{1-|x|^2}{2}x^i\phi_i=0,&
  \end{eqnarray*}
or equivalent to
   \begin{eqnarray}\nonumber\label{e5.3}
     &\xi^k=\xi^k(x,t),\ \ \eta=\eta(t), \ \ \phi=\alpha(x,t)u&\\ \nonumber
     &2\xi^i_{i}=-r\alpha(x,t)+\eta_t-\frac{4x^k\xi^k}{1-|x|^2}, \ \ \forall i&\\
     &\xi^i_j+\xi^j_i=0, \ \ \forall i\not=j&\\ \nonumber
     &\Big[(r+1-q)\alpha(x,t)-\eta_t\Big]u^q+\alpha_t(x,t)u^{r+1}-u\triangle_g\alpha=0&\\ \nonumber
     &u^r\partial_t\xi^k-\triangle_g\xi^k+\frac{(1-|x|^2)^2}{2}D_k\alpha-\frac{1-|x|^2}{2}\Big\{\eta_t+(1-n-r)\alpha\Big\}x^k&\\ \nonumber
     &+(2-n)\Bigg\{(x^i\xi^i)x^k-\frac{1-|x|^2}{2}\xi^k\Bigg\}=0, \ \ \forall k&
   \end{eqnarray}
As in Section 4, our results are divided into five cases:\\

\noindent\textbf{Case 1:} ($q, r+1, 1$ don't equal mutually)

By 4th-5th identities in \eqref{e5.3}, one has
   \begin{equation}\label{e5.4}
     \alpha=\frac{\beta}{r+1-q}, \ \ \eta(t)=\beta t+\gamma,\ \ \xi^k=\xi^k(x)
   \end{equation}
and
   \begin{equation}\label{e5.5}
     \begin{cases}
       \xi^i_i=-\frac{\beta(1-q)}{2(r+1-q)}-\frac{2x^k}{1-|x|^2}\xi^k, & \forall i\\
       \xi^i_j+\xi^j_i=0, & \forall i\not=j
     \end{cases}
   \end{equation}
and
   \begin{equation}\label{e5.6}
     \triangle_g\xi-(2-n)\frac{\xi}{1-z}=\Bigg\{\frac{\beta(2-n-q)}{r+1-q}+(2-n)y\cdot\xi\Bigg\}\frac{y}{(1-z)^2}, \ \ \xi\equiv(\xi^1,\cdots,\xi^n).
   \end{equation}
Applying Theorem \ref{t1.1} to \eqref{e5.5}, one concludes that
  \begin{equation}\label{e5.7}
    \alpha=0, \ \ \eta(t)=\gamma, \ \ \xi^i=\Sigma_{j\not=i}a^i_jx^j, \ \ a^i_j=-a^j_i, \ \ \forall i\not=j.
  \end{equation}
Therefore the infinite generators are spanned by
  \begin{eqnarray}\nonumber\label{e5.8}
    \overrightarrow{v}_1&=&\partial_t\\
    \overrightarrow{v}_2&=&\Sigma_{j\not=i}a^i_jx^j\partial_{x^i},
  \end{eqnarray}
and their corresponding group actions are generated by
   \begin{eqnarray}\nonumber\label{e5.9}
     g_1(\varepsilon) &: & (x,t,u)\to(\widetilde{x},\widetilde{t},\widetilde{u})=(x,t+\varepsilon, u)\\
     g_2(\varepsilon) &: & (x,t,u)\to(\widetilde{x},\widetilde{t},\widetilde{u})=(Ax,t, u),\ \ A\in SO(n)
   \end{eqnarray}
corresponding to translation in time and rotation on space.\\

\vspace{10pt}

\noindent\textbf{Case 2:} ($q=r+1\not=1$)

By 4th identity in \eqref{e5.3}, we get
   \begin{eqnarray}\nonumber\label{e5.10}
     &\alpha(x,t)=\eta(t)+\beta(x), \ \ \eta(t)=\kappa e^{rt}-\gamma,\ \ \xi^k=\xi^k(x)&\\
     &\triangle_g\beta=0&
   \end{eqnarray}
and
   \begin{equation}\label{e5.11}
     \begin{cases}
       \xi^i_i=\frac{r}{2}\Big[\gamma-\beta(x)\Big]-\frac{2x^k}{1-|x|^2}\xi^k, & \forall i\\
       \xi^i_j+\xi^j_i=0, & \forall i\not=j
     \end{cases}
   \end{equation}
and
   \begin{eqnarray}\nonumber\label{e5.12}
    &\triangle_g\xi-(2-n)\frac{\xi}{1-z}=\Bigg\{(1-n)\kappa e^{rt}+(1-n-r)\Big[\beta(x)-\gamma\Big]+(2-n)y\cdot\xi\Bigg\}\frac{y}{(1-z)^2}&\\
    &+\frac{2D\beta}{(1-z)^2},\ \ D\beta\equiv(\beta_1,\cdots,\beta_n).&
   \end{eqnarray}
Applying Theorem \ref{t1.1} to \eqref{e5.11}, one concludes that $\beta=\gamma$ and
   \begin{equation}\label{e5.13}
    \xi^k=a^k_jx^j+b^k, \ \ \forall k.
   \end{equation}
So, we get
  \begin{equation}\label{e5.14}
    \alpha=\kappa e^{rt}, \ \ \eta=\kappa e^{rt}-\gamma, \ \ \xi^k=\Sigma_{j\not=k}a^k_jx^j. \ \ a^i_j=-a^j_i, \ \ \forall i\not=j.
  \end{equation}
Thus, the infinite generators are spanned by
  \begin{eqnarray}\nonumber\label{e5.15}
    \overrightarrow{v}_1&=&e^{rt}\partial_t+e^{rt}u\partial_u\\
    \overrightarrow{v}_2&=&\partial_t\\ \nonumber
    \overrightarrow{v}_3&=&\Sigma_{j\not=i}a^i_jx^j\partial_{x^i},
  \end{eqnarray}
and their corresponding group actions are generated by
  \begin{equation}\label{e5.16}
     g_1(\varepsilon): \ \left(
    \begin{array}{c}
      x\\
      t \\
      u
    \end{array}
   \right)\to \left(
    \begin{array}{c}
      \widetilde{x}\\
      \widetilde{t}\\
      \widetilde{u}
    \end{array}
   \right)=   \left(
    \begin{array}{c}
      x\\
      -\frac{1}{r}\ln(e^{-rt}-r\varepsilon)\\
      u(e^{-rt}-r\varepsilon)^{-\frac{1}{r}}e^{-t}
    \end{array}
   \right)
  \end{equation}
and
   \begin{eqnarray}\nonumber\label{e5.17}
     g_2(\varepsilon) &: & (x,t,u)\to(\widetilde{x},\widetilde{t},\widetilde{u})=(x,t+\varepsilon, u)\\
     g_3(\varepsilon) &: & (x,t,u)\to(\widetilde{x},\widetilde{t},\widetilde{u})=(A_\varepsilon x, t, u), \ \ A_\varepsilon\in SO(n)
   \end{eqnarray}

\vspace{10pt}

\noindent\textbf{Case 3:} ($q=1\not=r+1$)

 We have $\xi^k=\xi^k(x)$ and
   \begin{eqnarray}\nonumber\label{e5.18}
      &\alpha=\alpha(x),\ \ \eta=\beta t+\gamma, \ \ \xi^k=\xi^k(x)&\\
      &\triangle_g\alpha-r\alpha+\beta=0&
   \end{eqnarray}
 and
   \begin{equation}\label{e5.19}
     \begin{cases}
       2\xi^i_i=-r\alpha(x)+\beta-\frac{4x^k\xi^k}{1-|x|^2}, & \forall i\\
       \xi^i_j+\xi^j_i=0, & \forall i\not=j
     \end{cases}
   \end{equation}
 and
    \begin{equation}\label{e5.20}
     \triangle_g\xi-(2-n)\frac{\xi}{1-z}=\Big\{\beta+(2-n-r)\alpha+y\cdot\xi\Big\}\frac{y}{(1-z)^2}+\frac{2D\alpha}{(1-z)^2}.
    \end{equation}
By Theorem \ref{t1.1} and \eqref{e5.19}, one obtains that
  \begin{equation}\label{e5.21}
    \alpha=\frac{\beta}{r},\ \ \eta=\beta t+\gamma, \ \ \xi^k=\Sigma_{j\not=k}a^k_jx^j, \ \ a^i_j=-a^j_i, \ \ \forall i\not=j.
  \end{equation}
Thus, the infinite generators are spanned by
  \begin{eqnarray}\nonumber\label{e5.22}
    \overrightarrow{v}_1&=&t\partial_t+u/r\partial_u\\
    \overrightarrow{v}_2&=&\partial_t\\ \nonumber
    \overrightarrow{v}_3&=&\Sigma_{j\not=i}a^i_jx^j\partial_{x^i},
  \end{eqnarray}
and their corresponding group actions are generated by
   \begin{eqnarray}\nonumber\label{e5.23}
     g_1(\varepsilon) &: & (x,t,u)\to(\widetilde{x},\widetilde{t},\widetilde{u})=(x,e^{\varepsilon}t, e^{\varepsilon/r}u)\\
     g_2(\varepsilon) &: & (x,t,u)\to(\widetilde{x},\widetilde{t},\widetilde{u})=(x,t+\varepsilon, u)\\ \nonumber
     g_3(\varepsilon) &: & (x,t,u)\to(\widetilde{x},\widetilde{t},\widetilde{u})=(A_\varepsilon x, t, u), \ \ A_\varepsilon\in SO(n)
   \end{eqnarray}

\vspace{10pt}

\noindent\textbf{Case 4:} ($r+1=1\not=q$)

One has $\xi^k=\xi^k(x,t)$ and
  \begin{eqnarray}\nonumber\label{e5.24}
    &&\begin{cases}
      (1-q)\alpha(x,t)=\eta_t\\
      \alpha_t-\frac{(1-|x|^2)^2}{4}\triangle\alpha-(n-2)\frac{1-|x|^2}{2}x\cdot\nabla\alpha=0
    \end{cases}\\
   &\Rightarrow& \eta=\beta t+\gamma, \ \ \alpha=\frac{\beta}{1-q}
  \end{eqnarray}
and
   \begin{equation}\label{e5.25}
     \begin{cases}
       2\xi^i_i=\beta-\frac{4x^k}{1-|x|^2}\xi^k, & \forall i\\
       \xi^i_j+\xi^j_i=0, & \forall i\not=j
     \end{cases}
   \end{equation}
and
  \begin{equation}\label{e5.26}
   \partial_t\xi-\triangle_g\xi+(2-n)\frac{\xi}{1-z}=-\Bigg\{\frac{\beta(2-n-q)}{1-q}+(2-n)y\cdot\xi\Bigg\}\frac{y}{(1-z)^2}.
  \end{equation}
Applying Theorem \ref{t1.1} to \eqref{e5.25}, one gets that
  \begin{equation}\label{e5.27}
    \alpha=0, \ \ \eta=\gamma, \ \ \xi^k=\Sigma_{j\not=k}a^k_j(t)x^j,\ \ a^k_j(t)=-a^j_k(t), \ \ \forall j\not=k.
  \end{equation}
Substituting into \eqref{e5.26}, one gets
   $$
    \Sigma_{j\not=k}\partial_ta^k_j(t)x^j=0, \ \ \forall k
   $$
or equivalently $\xi^k=\xi^k(x), \forall k$. Therefore the infinite generators are spanned by
  \begin{eqnarray}\nonumber\label{e5.28}
    \overrightarrow{v}_1&=&\partial_t\\
    \overrightarrow{v}_2&=&\Sigma_{j\not=i}a^i_jx^j\partial_{x^i},
  \end{eqnarray}
and their corresponding group actions are generated by
   \begin{eqnarray}\nonumber\label{e5.29}
     g_1(\varepsilon) &: & (x,t,u)\to(\widetilde{x},\widetilde{t},\widetilde{u})=(x,t+\varepsilon, u)\\
     g_2(\varepsilon) &: & (x,t,u)\to(\widetilde{x},\widetilde{t},\widetilde{u})=(Ax,t, u),\ \ A\in SO(n)
   \end{eqnarray}
corresponding to translation in time and rotation on space respectively.\\

\vspace{10pt}

\noindent\textbf{Case 5:} ($q=r+1=1$)

We have
 \begin{equation}\label{e5.30}
   \alpha_t-\triangle_g\alpha=\eta_t
 \end{equation}
and
   \begin{equation}\label{e5.31}
     \begin{cases}
       2\xi^i_i=\eta_t-\frac{4x^k}{1-|x|^2}\xi^k, & \forall i\\
       \xi^i_j+\xi^j_i=0, & \forall i\not=j
     \end{cases}
   \end{equation}
and
  \begin{equation}\label{e5.32}
  \partial_t\xi-\triangle_g\xi+(2-n)\frac{\xi}{1-z}=-\Bigg\{\eta_t+(1-n)\alpha+(2-n)y\cdot\xi\Bigg\}\frac{y}{(1-z)^2}-\frac{2D\alpha}{(1-z)^2}.
  \end{equation}
By Theorem \ref{t1.1}, it's inferred from \eqref{e5.30}-\eqref{e5.31} that
  \begin{equation}\label{e5.33}
    \alpha_t-\triangle_g\alpha=0, \ \ \eta=\gamma,\ \ \xi^k=\Sigma_{j\not=k}a^k_j(t)x^j, \ \ a^k_j(t)=-a^j_k(t), \ \  \forall j\not=k.
  \end{equation}
Combining with \eqref{e5.32}, one concludes that
  \begin{equation}\label{e5.34}
   D_k\alpha+\frac{(n-1)x^k}{1-|x|^2}\alpha=-\Sigma_{j\not=k}\frac{(1-|x|^2)^2}{2}\partial_ta^k_j(t)x^j.
  \end{equation}
Taking a second derivative on $x^l$, there holds
  \begin{eqnarray}\nonumber\label{e5.35}
    D_{kl}\alpha&=&-(n-1)\Bigg\{\frac{\delta_{kl}}{1-|x|^2}+\frac{2x^kx^l}{(1-|x|^2)^2}\Bigg\}\alpha-(n-1)\frac{x^k}{1-|x|^2}D_l\alpha\\
     &&+2(1-|x|^2)\Sigma_{j\not=k, j\not=l}x^lx^j\partial_ta^k_j(t)-\Bigg\{\frac{(1-|x|^2)^2}{2}-2(x^l)^2(1-|x|^2)\Bigg\}\partial_ta^k_l(t).
  \end{eqnarray}
Substituting \eqref{e5.34} into \eqref{e5.35}, we derive that
   \begin{eqnarray*}
    &D_{kl}\alpha+(n-1)\Bigg\{\frac{\delta_{kl}}{1-|x|^2}-(n-3)\frac{x^kx^l}{(1-|x|^2)^2}\Bigg\}\alpha&\\
    &=\frac{n-1}{2}(1-|x|^2)\Sigma_{j\not=k, j\not=l}x^kx^j\partial_ta^l_j(t)+2(1-|x|^2)\Sigma_{j\not=k, j\not=l}x^lx^j\partial_ta^k_j(t)&\\
    &-\frac{n-1}{2}(1-|x|^2)(x^k)^2\partial_ta^k_k(t)-(1-|x|^2)\Bigg\{\frac{1-|x|^2}{2}-2(x^l)^2\Bigg\}\partial_ta^k_l(t).&
   \end{eqnarray*}
Exchanging the indices $k,l$ and using the symmetry of second derivatives, it yields that
   \begin{eqnarray*}
     &\frac{5-n}{2}\Sigma_{j\not=k, j\not=l}x^j\Bigg\{x^k\partial_ta^l_j(t)-x^l\partial_ta^k_j(t)\Bigg\}&\\
     &+\Bigg\{(1-|x|^2)+\frac{5-n}{2}\Big[(x^k)^2+(x^l)^2\Big]\Bigg\}\partial_ta^k_l(t)=0, \ \ \forall k\not=l&.
   \end{eqnarray*}
Consequently, we have
   \begin{equation}\label{e5.36}
     a^k_j(t)=a^k_j, \ \ \forall j\not=k
   \end{equation}
and so
   \begin{equation}\label{e5.37}
     \alpha=\beta, \ \ \eta=\gamma, \xi^k=\Sigma_{j\not=k}a^k_jx^j, \ \ a^k_j=-a^j_k, \ \ \forall j\not=k.
   \end{equation}
Thus, the infinite generators are spanned by
  \begin{eqnarray}\nonumber\label{e5.38}
    \overrightarrow{v}_1&=&u\partial_u\\
    \overrightarrow{v}_2&=&\partial_t\\ \nonumber
    \overrightarrow{v}_3&=&\Sigma_{j\not=i}a^i_jx^j\partial_{x^i},
  \end{eqnarray}
and their corresponding group actions are generated by
   \begin{eqnarray}\nonumber\label{e5.39}
     g_1(\varepsilon) &: & (x,t,u)\to(\widetilde{x},\widetilde{t},\widetilde{u})=(x,t,e^{\varepsilon}u)\\
     g_2(\varepsilon) &: & (x,t,u)\to(\widetilde{x},\widetilde{t},\widetilde{u})=(x,t+\varepsilon,u)\\ \nonumber
     g_3(\varepsilon) &: & (x,t,u)\to(\widetilde{x},\widetilde{t},\widetilde{u})=(A_\varepsilon x, t, u), \ \ A_\varepsilon\in SO(n)
   \end{eqnarray}

\vspace{40pt}

\section*{Acknowledgments}
The author SZ would like to express his deepest gratitude to Professors Kai-Seng Chou and Xu-Jia Wang for their constant encouragements and warm-hearted helps. Special thanks were also owed to Professor Cheng-Jie Yu for valuable conversation.

\end{document}